\documentclass[preprint,12pt]{elsarticle}
\usepackage{amsmath,amssymb,amsfonts}
\usepackage{amsthm}

\newtheorem{theorem}{Theorem}

\newtheorem{lemma}[theorem]{Lemma}
\newtheorem{definition}{Definition}
	
\journal{Journal of Combinatorial Theory, Series A}

\begin{document}

\begin{frontmatter}



\title{MacMahon-type $q$-series}


\author{Mircea Merca}

\affiliation{organization={Department of Mathematical Methods and Models, Fundamental Sciences Applied in Engineering Research Center, National University of Science and Technology Politehnica Bucharest},%Department and Organization
            %addressline={}, 
            city={Bucharest},
            postcode={060042}, 
            %state={},
            country={Romania}}

\affiliation{organization={Academy of Romanian Scientists},%Department and Organization
	%addressline={}, 
	city={Bucharest},
	postcode={050044}, 
	%state={},
	country={Romania}}
	
\begin{abstract}
Motivated by earlier work of P.~A.~MacMahon and recent contributions of T.~Amdeberhan, G.~E.~Andrews, K.~Ono, A.~Singh, and R.~Tauraso on higher-order partition enumerants, we study a class of $q$-series arising from nested divisor 
structures. In particular, we consider the $q$-series
\[
V_k(q)
= \sum_{1 \le n_1 \le n_2 \le \cdots \le n_k}
\frac{q^{\,n_1+n_2+\cdots+n_k}}
{(1-q^{n_1})^2(1-q^{n_2})^2\cdots(1-q^{n_k})^2},
\]
introduced recently as MacMahon-type generating functions. We further define a new MacMahon-type series
\[
W_k(q)
= \sum_{1 \le n_1 \le n_2 \le \cdots \le n_k}
\frac{q^{\,2(n_1+n_2+\cdots+n_k)-k}}
{(1-q^{2n_1-1})^2(1-q^{2n_2-1})^2\cdots(1-q^{2n_k-1})^2},
\]
and establish families of identities, generating function relations, and hypergeometric 
representations for the truncated forms of $V_k(q)$ and $W_k(q)$. Connections with overpartition pairs and bipartitions with distinct odd parts arise naturally in this context.
\end{abstract}



\begin{keyword}
	partitions \sep overpartitions \sep $q$-series



\MSC[2010] 11P81 \sep 11P82 \sep 05A19 \sep 05A20
\end{keyword}

\end{frontmatter}



\section{Introduction}
\label{S1}

The partitions of an integer represent a fascinating subject at the intersection of combinatorics and number theory. A partition of a positive integer $n$  is a way of writing  $n$ as a sum of positive integers, where the order of the terms (called parts) does not matter \cite{Andrews98}.
The multiplicity of parts provides a useful alternative representation for the partitions of a positive integer $n$:
$$
n = t_1 + 2t_2 + 3t_3 + \cdots + nt_n,
$$
where $t_i$ is the number of times the part $i$ appears. These counts $t_i$ are non-negative integers (possibly zero), and this form is often encoded compactly as
$$
(1^{t_1} 2^{t_2} \ldots n^{t_n})\qquad\text{or}\qquad (n^{t_n}\ldots 2^{t_2}1^{t_1}),
$$
which explicitly shows the multiplicity of each part. 
 For instance, the integer $4$ has five partitions, represented as:
$$
(4^1), \quad (3^1 1^1), \quad (2^2), \quad (2^1 1^2), \quad (1^4).
$$
Note that in this frequency-based description, it is possible that all $t_i = 0$; 
this case corresponds to the \textit{empty partition} of $n = 0$, conventionally denoted by the empty list $()$. 
Hence, the definitions are consistent when we regard $()$ as the unique partition of zero.

\smallskip
A partition of an integer $n$ is said to involve parts of $k$ distinct magnitudes if the set of values taken by its parts has a cardinality exactly equal to $k$. For example, if we consider partitions of $8$ involving parts of $3$ different magnitudes, we might have:
\begin{align}
	(5^1 2^1 1^1),\ (4^1 3^1 1^1),\ (4^1 2^1 1^2),\ (3^1 2^2 1^1),\ (3^1 2^1 1^3).\label{eq2}    
\end{align}
The enumeration and structural analysis of these partitions has established deep connections between combinatorial counting problems and sophisticated analytical tools, particularly $q$-series and modular forms. 
Considering partitions involve parts of $k$ distinct magnitudes, MacMahon introduced 
two fundamental families of partition functions:
$$
 		a^{\pm}_{k}(n) := \sum_{\substack{\lambda_1t_1+\lambda_2t_2+\cdots +\lambda_kt_k = n \\ 1\leq\lambda_1<\lambda_2<\cdots<\lambda_k\\ (t_1,t_2,\ldots,t_k)\in\mathbb{N}^k}} (\pm 1)^{t_1+t_2+\cdots+t_k+k}\, t_1t_2\cdots t_k
$$
and
$$
 		c^{\pm}_{k}(n) := \sum_{\substack{(2\lambda_1-1)t_1+(2\lambda_2-1)t_2+\cdots +(2\lambda_k-1)t_k = n \\ 1\leq\lambda_1<\lambda_2<\cdots<\lambda_k\\ (t_1,t_2,\ldots,t_k)\in\mathbb{N}^k}} (\pm 1)^{t_1+t_2+\cdots+t_k+k}\, t_1t_2\cdots t_k.
$$
These serve as the basis for two significant families of $q$-series, which MacMahon \cite{MacMahon} established as their generating functions:
$$
A^{\pm}_k(q) = \sum_{n=0}^\infty a^{\pm}_{k}(n)\,q^n:=\sum_{0<n_1<n_2<\cdots <n_k}
\prod_{i=1}^k \frac{q^{n_i}}{(1\mp q^{n_i})^2}
$$
and    
$$
C^{\pm}_k(q) = \sum_{n=0}^\infty c^{\pm}_{k}(n)\,q^n:=\sum_{0<n_1<n_2<\cdots <n_k}
\prod_{i=1}^k \frac{q^{2n_i-1}}{(1\mp q^{2n_i-1})^2}.
$$
We remark that $A^{\pm}_0(q)$ and $C^{\pm}_0(q)$ are defined to be $1$.

\smallskip
These functions have recently experienced a resurgence of interest, resulting in a considerable number of new articles that highlight their continued mathematical importance \cite{Amdeberhan23,Amdeberhan24,Amdeberhan25,Andrews13,Andrews22,Bachmann,Liu25,Merca2025,Ono,Rose,Sellers25,Xia25}.
Amdeberhan et al. \cite{Amdeberhan24} recently refined MacMahon's original scope, suggesting that while the established approach is sound, the study could be expanded to include partitions with exactly $k$ size of parts, possibly equal. Following this idea, we consider two families of $q$-series:
$$
V^{\pm}_k(q) = \sum_{n=0}^\infty v^{\pm}_k(n)\,q^n:=\sum_{1\leq n_1\leq n_2\leq \cdots \leq n_k}
\prod_{i=1}^k \frac{q^{n_i}}{(1\mp q^{n_i})^2},
$$
and
$$
W^{\pm}_k(q) = \sum_{n=0}^\infty w^{\pm}_k(n)\,q^n:=\sum_{1\leq n_1\leq n_2\leq \cdots \leq n_k}
\prod_{i=1}^k \frac{q^{2n_i-1}}{(1\mp q^{2n_i-1})^2},
$$
where the families of partition functions $v^{\pm}_k(n)$ and $w^{\pm}_k(n)$ are defined by
$$
v^{\pm}_k(n):= \sum_{\substack{\lambda_1t_1+\lambda_2t_2+\cdots +\lambda_kt_k = n \\ 1\leq\lambda_1\leq\lambda_2\leq\cdots\leq\lambda_k\\ (t_1,t_2,\ldots,t_k)\in\mathbb{N}^k}} (\pm 1)^{t_1+t_2+\cdots+t_k+k}\, t_1t_2\cdots t_k
$$
and
$$
w^{\pm}_k(n):= \sum_{\substack{(2\lambda_1-1)t_1+(2\lambda_2-1)t_2+\cdots +(2\lambda_k-1)t_k = n \\ 1\leq\lambda_1\leq\lambda_2\leq\cdots\leq\lambda_k\\ (t_1,t_2,\ldots,t_k)\in\mathbb{N}^k}} (\pm 1)^{t_1+t_2+\cdots+t_k+k}\, t_1t_2\cdots t_k.
$$
For example, if we consider partitions of $5$ involving parts of $3$ magnitudes (possibly equal), we might have:
\begin{align*}
	(3^1 1^1 1^1),\ (2^1 2^1 1^1),\ (2^1 1^1 1^2),\ (2^1 1^2 1^1),\ (1^1 1^1 1^3),\\
	(1^1 1^2 1^2),\ (1^1 1^3 1^1),\ (1^2 1^1 1^2),\ (1^2 1^2 1^1),\ (1^3 1^1 1^1).  
\end{align*}
So we have:
\begin{align*}
	v^+_3(5) = &1\cdot 1\cdot 1+ 1\cdot 1\cdot 1+1\cdot 1\cdot 2+ 1\cdot 2\cdot 1+1\cdot 1\cdot 3\\
	&+ 1\cdot 2\cdot 2+1\cdot 3\cdot 1+2\cdot 1\cdot 2+2\cdot 2\cdot 1+3\cdot 1\cdot 1 = 27\\
	v^-_3(5) = &1\cdot 1\cdot 1+ 1\cdot 1\cdot 1-1\cdot 1\cdot 2- 1\cdot 2\cdot 1+1\cdot 1\cdot 3\\
	&+ 1\cdot 2\cdot 2+1\cdot 3\cdot 1+2\cdot 1\cdot 2+2\cdot 2\cdot 1+3\cdot 1\cdot 1 = 19.
\end{align*}
On the other hand, partitions of $5$ involving odd parts of $3$ magnitudes (possibly equal) are:
\begin{align*}
	(3^1 1^1 1^1),\ (1^1 1^1 1^3),\ 
	(1^1 1^2 1^2),\ (1^1 1^3 1^1),\ (1^2 1^1 1^2),\ (1^2 1^2 1^1),\ (1^3 1^1 1^1).  
\end{align*}
So we have:
\begin{align*}
	w^{+}_{3}(5) = &1\cdot 1\cdot 1+ 1\cdot 1\cdot 3+ 1\cdot 2\cdot 2\\
	&+1\cdot 3\cdot 1+2\cdot 1\cdot 2+2\cdot 2\cdot 1+3\cdot 1\cdot 1 = 22.
\end{align*}
We remark that
$$
w^-_k(n)=(-1)^{n+k}\,w^+_k(n).
$$

In this paper, in order to explore properties of the families of MacMahon-type $q$-series $V_k^{\pm}(q)$ and $W_k^{\pm}(q)$, we introduce the following partitions functions.

\begin{definition}
	For positive integers $k$, $m$ and $n$, we define:
\begin{enumerate}
	\item[(i)] $\displaystyle{
v^{\pm}_{k,m}(n):= \sum_{\substack{\lambda_1t_1+\lambda_2t_2+\cdots +\lambda_kt_k = n \\ 1\leq\lambda_1\leq\lambda_2\leq\cdots\leq\lambda_k\leq m\\ (t_1,t_2,\ldots,t_k)\in\mathbb{N}^k}} (\pm 1)^{t_1+t_2+\cdots+t_k+k}\, t_1t_2\cdots t_k;}$
	\item[(ii)] $\displaystyle{
w^{\pm}_{k,m}(n):= \sum_{\substack{(2\lambda_1-1)t_1+(2\lambda_2-1)t_2+\cdots +(2\lambda_k-1)t_k = n \\ 1\leq\lambda_1\leq\lambda_2\leq\cdots\leq\lambda_k\leq m\\ (t_1,t_2,\ldots,t_k)\in\mathbb{N}^k}} (\pm 1)^{t_1+t_2+\cdots+t_k+k}\, t_1t_2\cdots t_k.}$
\end{enumerate}
\end{definition}

For positive integers $k$ and $m$, we remark that the generating functions of $v^{\pm}_{k,m}(n)$ and $w^{\pm}_{k,m}(n)$ are truncated forms of the MacMahon-type $q$-series $V^{\pm}_{k}(q)$ and $W^{\pm}_{k}(q)$, i.e.,
$$
V^{\pm}_{k,m}(q) = \sum_{n=0}^\infty v^{\pm}_{k,m}(n)\,q^n:=\sum_{1\leq n_1\leq n_2\leq \cdots \leq n_k\leq m}
\prod_{i=1}^k \frac{q^{n_i}}{(1\mp q^{n_i})^2},
$$
and
$$
W^{\pm}_{k,m}(q) = \sum_{n=0}^\infty w^{\pm}_{k,m}(n)\,q^n:=\sum_{1\leq n_1\leq n_2\leq \cdots \leq n_k \leq m}
\prod_{i=1}^k \frac{q^{2n_i-1}}{(1\mp q^{2n_i-1})^2},
$$
with $V^{\pm}_{0,m}(q)$ and $W^{\pm}_{0,m}(q)$ defined to be $1$.

\medskip
In this context, we recall the definition of the $q$-binomial coefficient, also known as Gaussian polynomial. For integers 
$m,k \in \mathbb{Z}$, it is defined by
\[
{m \brack k} = {m \brack k}_q :=
\begin{cases}
\dfrac{(q;q)_m}{(q;q)_k\, (q;q)_{m-k}}, & \text{if } 0 \le k \le m,\\[6pt]
0, & \text{otherwise}.
\end{cases}
\]
Here $(a;q)_n$ denotes the $q$-shifted factorial (the $q$-Pochhammer symbol), 
defined for $n \in \mathbb{Z}_{\ge 0}$ by
\[
(a;q)_n = \prod_{r=0}^{n-1} (1 - a q^{\,r}), \qquad (a;q)_0 = 1.
\]
We also use the infinite $q$-product
\[
(a;q)_\infty = \lim_{n\to\infty} (a;q)_n,
\]
which exists and converges whenever $|q| < 1$.

\medskip
The first result of this paper exhibits two infinite families of linear combinations of the $q$-series
$V_{k,m}^{\pm}(q)$ and $W_{k,m}^{\pm}(q)$, expressed in terms of Gaussian polynomials or basic hypergeometric series \cite{gasper}. To this end, we define
\begin{align*}
H_{k,m}(q,z) &:=\sum_{j=0}^\infty {m-1+j \brack j}_q {m-1+k+j \brack k+j}_q z^{j} \\
&= \frac{(q^m;q)_k}{(q;q)_k}\, {_{2}}\phi_{1}\bigg(\begin{matrix}q^m, q^{m+k}\\ q^{k+1}\end{matrix}\,;q,z \bigg),
\end{align*}
where
$$
{_{2}}\phi_{1}\bigg(\begin{matrix}a, b\\ c\end{matrix}\,;q,z\bigg) 
=\sum_{n=0}^\infty \frac{(a;q)_n(b;q)_n}{(q;q)_n(c;q)_n} z^n.
$$

\medskip
\begin{theorem}\label{T1}
	Let $k$ and $m$ be non-negative integers. For $|q|<1$,
	\begin{enumerate}
		\item [(i)] $\displaystyle{\sum_{j=k}^\infty (\mp 1)^{j-k} \binom{2j}{j-k} V^{\pm}_{j,m}(q) = (\pm q;q)^2_m\,q^k\, H_{k,m}(q,q^2)}$
	%	\item [] $\displaystyle{\qquad =(\pm q;q)^2_m \sum_{j=0}^\infty 
	%		{m-1+j \brack j} {m-1+k+j \brack k+j} q^{2j+k}}$;
		\item [(ii)] $\displaystyle{\sum_{j=k}^\infty (\mp 1)^{j-k} \binom{2j}{j-k} W^{\pm}_{j,m}(q)=(\pm q;q^2)^2_m\,q^k\,H_{k,m}(q^2,q^2)}.$
	%	\item [] $\displaystyle{\qquad = (\pm q;q^2)^2_m \sum_{j=0}^\infty 
	%		{m-1+j \brack j}_{q^2} {m-1+k+j \brack k+j}_{q^2} q^{2j+k}}$.		
	\end{enumerate}
\end{theorem}
\medskip

By Theorem \ref{T1},  we deduce that the $q$-series
$V^{\pm}_{j,m}(q)$ and $W^{\pm}_{j,m}(q)$ admit explicit expressions in terms of $q$-binomial coefficients.
To formulate these expressions, we introduce the inverse coefficients $B_{k,j}$, defined by
\[
B_{0,0}=1, \qquad B_{k,0}=2 \quad (k \geq 1),
\]
and, for $j \geq 1$ and $k \geq j$,
\[
B_{k,j} = \frac{2k}{k+j} \binom{k+j}{2j}.
\]
With these coefficients, the $q$-series can be written in closed form, providing a direct combinatorial interpretation.

\medskip
\begin{theorem}\label{TT2}
	Let $j$ and $m$ be non-negative integers. For $|q|<1$,
	\begin{enumerate}
		\item [(i)] $\displaystyle{ V^{\pm}_{j,m}(q) 
		= (\pm q;q)_m^2 \sum_{k=j}^\infty (\pm 1)^{k-j} B_{k,j}\,q^k\,H_{k,m}(q,q^2)}$,\\
		% \sum_{\ell=0}^{\infty} {m-1+\ell \brack \ell}{m-1+k+\ell \brack k+\ell} q^{2\ell+k},\\
		\item [(ii)] $\displaystyle{ W^{\pm}_{j,m}(q) 
	= (\pm q;q^2)^2_m \sum_{k=j}^\infty (\pm 1)^{k-j} B_{k,j}\,q^k\,H_{k,m}(q^2,q^2)}$. 
	%\sum_{\ell=0}^{\infty} {m-1+\ell \brack \ell}_{q^2}{m-1+k+\ell \brack k+\ell}_{q^2} q^{2\ell+k}.
\end{enumerate}
\end{theorem}
\medskip

In light of Theorems~\ref{T1} and~\ref{TT2}, we are now able to deduce structural properties of the MacMahon-type $q$-series 
$V^{\pm}_k(q)$ and $W^{\pm}_k(q)$.
In particular, taking the limit $m\to\infty$ in Theorem~\ref{T1} yields the following result.

\medskip
\begin{theorem}\label{T4}
	Let $k$ be non-negative integer. For $|q|<1$,
	\begin{enumerate}
		\item [(i)] $\displaystyle{\sum_{j=k}^\infty (\mp 1)^{j-k} \binom{2j}{j-k} V^{\pm}_{j}(q)}$
		\item [] $\displaystyle{\qquad\qquad = \frac{(\pm q;q)^2_\infty}{(q;q)^2_\infty} \left(-1+(1+q^k) \sum_{j=k}^\infty (-1)^{j-k}\,q^{j(j+1)/2-k(k+1)/2}\right)}$;
		\item [(ii)] $\displaystyle{\sum_{j=k}^\infty (\mp 1)^{j-k} \binom{2j}{j-k} W^{\pm}_{j}(q) = \frac{(\pm q;q^2)^2_\infty}{(q^2;q^2)^2_\infty}\sum_{j=k}^\infty (-1)^{j-k}\,q^{j(j+1)-k^2}.}$
	\end{enumerate}
\end{theorem}
\medskip

%Recall that, an \emph{overpartition} of \(n\) is a partition of \(n\) in which the first occurrence 
%of a number may be overlined. For example, there are \(14\) overpartitions of \(4\):
%\begin{align*}
%& 4,\ \overline{4},\ 3+1,\ \overline{3}+1,\ 3+\overline{1},\ \overline{3}+\overline{1},\
%2+2,\ \overline{2}+2,\ 2+\overline{2},\\
%& \overline{2}+\overline{2},\
%2+1+1,\ \overline{2}+1+1,\ 2+\overline{1}+1,\ \overline{2}+\overline{1}+1.
%\end{align*}
%An \emph{overpartition pair} of \(n\) is an ordered pair \((\mu,\lambda)\) of 
%overpartitions such that the total sum of all parts in both \(\mu\) and \(\lambda\) 
%is \(n\). For example, there are \(12\) overpartition pairs of \(2\):
%\[
%(2,\emptyset),\ (\overline{2},\emptyset),\ 
%(1+1,\emptyset),\ (\overline{1}+1,\emptyset),\ 
%(1,1),\ (\overline{1},1),\ (1,\overline{1}),\ (\overline{1},\overline{1}),
%\]
%\[
%(\emptyset,2),\ (\emptyset,\overline{2}),\ 
%(\emptyset,1+1),\ (\emptyset,\overline{1}+1).
%\]
These equalities separate a product part and a theta/false theta part given by alternating sums of triangular powers.  We known that
\[
\frac{(-q;q)_\infty^2}{(q;q)_\infty^2}
= \sum_{n=0}^\infty \overline{pp}(n)\,q^n,
\]
where \(\overline{pp}(n)\) denotes the number of overpartition pairs of \(n\)
\cite{Lovejoy}. Likewise,
\[
\frac{(-q;q^2)_\infty^2}{(q^2;q^2)_\infty^2}
= \sum_{n=0}^\infty pod_{-2}(n)\,q^n,
\]
where \(pod_{-2}(n)\) denotes the number of bipartitions of \(n\) with odd parts
distinct and even parts are unrestricted \cite{Chen}. By the proof of Theorem \ref{T4}, we see that the coefficients of $q^n$ in the $q$-series
$$\frac{(-q;q)^2_\infty}{(q;q)^2_\infty} \left(-1+(1+q^k) \sum_{j=k}^\infty (-1)^{j-k}\,q^{j(j+1)/2-k(k+1)/2}\right)$$
and
$$\frac{(-q;q^2)^2_\infty}{(q^2;q^2)^2_\infty}\sum_{j=k}^\infty (-1)^{j-k}\,q^{j(j+1)-k^2}$$
are all non-negative. Thus we deduce that, for the non-negative integers $k$ and $n$, the partition inequality
$$
\sum_{j=k}^\infty \binom{2j}{j-k}v^{-}_j(n)\geq 0
$$
is equivalent to
$$
-\overline{pp}(n)+\sum_{j=k}^\infty (-1)^{j-k} \Big(\overline{pp}(n-T_j+T_k)+\overline{pp}(n-T_j+T_{k-1}) \Big)\geq 0,
$$
where $T_n=n(n+1)/2$ is the $n$-th triangular number. On the other hand, for the non-negative integers $k$ and $n$, the partition inequality
$$
\sum_{j=k}^\infty \binom{2j}{j-k}w^{-}_j(n)\geq 0
$$
is equivalent to
$$
\sum_{j=k}^\infty (-1)^{j-k}\,pod_{-2}(n-T_j+k^2) \geq 0.
$$
In the specialisation $k=0$, Theorem \ref{T4} collapses to the pair of identities
\begin{align*}
	& \sum_{j=0}^\infty (\mp 1)^{j} \binom{2j}{j} V^{\pm}_{j}(q) 
	= \frac{(\pm q;q)^2_\infty}{(q;q)^2_\infty} \left(-1+2 \sum_{j=0}^\infty (-1)^{j}\,q^{j(j+1)/2}\right),\\
	& \sum_{j=0}^\infty (\mp 1)^{j} \binom{2j}{j} W^{\pm}_{j}(q) = \frac{(\pm q;q^2)^2_\infty}{(q^2;q^2)^2_\infty}\sum_{j=0}^\infty (-1)^{j}\,q^{j(j+1)}.
\end{align*}
In this context, it is natural to ask: what combinatorial object is enumerated by the coefficient of $q^n$ in
\[
\frac{(-q;q)^2_\infty}{(q;q)^2_\infty} \Biggl(-1 + 2 \sum_{j=0}^\infty (-1)^{j}\,q^{j(j+1)/2}\Biggr)?
\]
Similarly, what does the coefficient of $q^n$ in
\[
\frac{(-q;q^2)^2_\infty}{(q^2;q^2)^2_\infty} \sum_{j=0}^\infty (-1)^{j}\,q^{j(j+1)}
\]
count?

In the limiting case $m \to \infty$ of Theorem~\ref{TT2}, the partition functions $v_k^{-}(n)$ and $w_k^{\pm}(n)$ admit interpretations, respectively, in terms of overpartition pairs and bipartitions with odd parts.

\medskip
\begin{theorem}\label{TT4}
	Let $j$ be a non-negative integer. For $|q|<1$,
%	\begin{align*}
%		& V^{\pm}_{j}(q) 
%		= \frac{(\pm q;q)^2_\infty}{(q;q)^2_\infty} \sum_{k=j}^\infty (\pm 1)^{k-j} B_{k,j} 
%		\Big(-1+(1+q^k) \sum_{\ell=k}^\infty (-1)^{\ell-k}\,q^{\ell(\ell+1)/2-k(k+1)/2}\Big),\\
%		& W^{\pm}_{j}(q) 
%		= \frac{(\pm q;q^2)^2_\infty}{(q^2;q^2)^2_\infty} \sum_{k=j}^\infty (\pm 1)^{k-j} B_{k,j} 
%		\sum_{\ell=k}^\infty (-1)^{\ell-k}\,q^{\ell(\ell+1)-k^2}.
%	\end{align*}
	\begin{enumerate}
		\item [(i)] $\displaystyle{V^{\pm}_{j}(q) 
			= \frac{(\pm q;q)^2_\infty}{(q;q)^2_\infty} \times}$
		\item [] $\displaystyle{ \qquad \times \sum_{k=j}^\infty (\pm 1)^{k-j} B_{k,j} 
			\Big(-1+(1+q^k) \sum_{\ell=k}^\infty (-1)^{\ell-k}\,q^{\ell(\ell+1)/2-k(k+1)/2}\Big);}$
		\item [(ii)] $\displaystyle{W^{\pm}_{j}(q) 
			= \frac{(\pm q;q^2)^2_\infty}{(q^2;q^2)^2_\infty} \sum_{k=j}^\infty (\pm 1)^{k-j} B_{k,j} 
			\sum_{\ell=k}^\infty (-1)^{\ell-k}\,q^{\ell(\ell+1)-k^2}.}$
	\end{enumerate}
\end{theorem}
\medskip

In the special case $j=0$, the theorem yields two identities that express
products of Ramanujan theta functions \cite[Theorem 3.5]{Cooper}
$$
\varphi(-q)=\frac{(q;q)_\infty}{(-q;q)_\infty}\qquad\text{and}\qquad\psi(q)=\frac{(q^2;q^2)_\infty}{(q;q^2)_\infty}
$$
in terms of the coefficients $B_{k,0}$:
\begin{align*}
\Big( \sum_{k=-\infty}^{\infty} (-q)^{k^2} \Big)^{\!2}
&= \sum_{k=0}^{\infty} B_{k,0}\,\Big( (-1)^{k+1}
+ (1+q^k)\sum_{j=k}^{\infty} (-1)^j\,
q^{\,j(j+1)/2-k(k+1)/2} \Big), \\
\Big( \sum_{k=0}^{\infty} q^{k(k+1)/2} \Big)^{\!2}
&= \sum_{k=0}^{\infty} B_{k,0}
\sum_{j=k}^{\infty} (-1)^{\,j-k}\,
q^{\,j(j+1)-k^2}.
\end{align*}
Since each left–hand side is a square of a Ramanujan theta function,
these identities give a decomposition of their coefficients into
alternating sums depending on the parameters $B_{k,0}$.
In particular, the coefficient of $q^n$ on the left--hand side counts
the number of ordered representations of $n$
as the sum of two squares with a prescribed sign pattern or, respectively, 
the sum of two triangular numbers.
Thus, the case $j=0$ of theorem provides new weighted expansions for the
representation numbers encoded by these classical theta functions.

The remainder of the paper is devoted to the proofs of Theorems \ref{T1}–\ref{T4}.
Theorem \ref{TT4} follows immediately from Theorems \ref{TT2} and \ref{T4}.

\section{Proof of Theorem \ref{T1}}

\allowdisplaybreaks{
	The following result was proved by Cauchy (\cite[Theorem 26]{Johnson}): if $n$ is any nonnegative integer and $|q|$ and $|t|$ are both less than $1$, then
	\begin{align}
		\sum_{k=0}^\infty {n-1+k \brack k} t^k = \frac{1}{(t;q)_n}.\label{Cauchy}
	\end{align}
	Considering this identity, we can write:
	\begin{align*}
		& \sum_{i=0}^\infty V^{\pm}_{i,m}(q)\, (z\mp z^{-1})^{2i} \\ 
		&\qquad   = \prod_{i=1}^m \left(1-\frac{q^i}{(1\mp q^i)^2} (z\mp z^{-1})^2 \right)^{-1}\\
		%&\qquad   = (\pm q;q)^2_m \prod_{i=1}^m (1-q^iz^2)^{-1}(1-q^iz^{-2})^{-1} \\
		&\qquad   = (\pm q;q)^2_m \prod_{i=1}^m \frac{1}{(1-q^iz^2)(1-q^iz^{-2})} \\
		&\qquad   = \frac{(\pm q;q)^2_m}{(qz^2;q)_m\,(qz^{-2};q)_m} \\
		%&\qquad   = \frac{(qz^{-2};q)_m}{(\pm q;q)^2_m} \sum_{i=0}^m (-1)^i\, q^{\binom{i+1}{2}} {m \brack i} z^{2i} \\
		%\intertext{\hfill by \eqref{Rbt} with $t$ replaced by $-qz^2$}
		&\qquad   = (\pm q;q)^2_m \left(\sum_{i=0}^\infty  {m-1+i \brack i} q^i\, z^{2i} \right)
		\left(\sum_{i=0}^\infty  {m-1+i \brack i} q^i\,z^{-2i} \right)\\
		\intertext{\hfill by \eqref{Cauchy} with $t$ replaced by $qz^{\pm2}$}
		&\qquad   = (\pm q;q)^2_m \sum_{j=0}^\infty {m-1+j \brack j}^2 q^{2j}\\ 
		&\qquad\quad +(\pm q;q)^2_m \sum_{i=1}^\infty (z^{2i}+z^{-2i}) \sum_{j=0}^\infty {m-1+j \brack j} {m-1+i+j \brack i+j} q^{i+2j}\\
		&\qquad = (\pm q;q)^2_m \sum_{i=-\infty}^\infty z^{2i} \sum_{j=0}^\infty 
		{m-1+j \brack j} {m-1+|i|+j \brack |i|+j} q^{|i|+2j}.
	\end{align*}
	On the other hand, we have:
	\begin{align*}
		& \sum_{i=0}^\infty V^{\pm}_{i,m}(q)\, (z\mp z^{-1})^{2i} \\ 
		& \qquad = \sum_{i=0}^\infty \sum_{j=0}^{2i} (\mp 1)^j \,V^{\pm}_{i,m}(q) \binom{2i}{j} z^{2i-2j}\\
		& \qquad = \sum_{j=0}^\infty (\mp 1)^j \binom{2j}{j} V^{\pm}_{j,m}(q) + 
		\sum_{i=1}^\infty (z^{2i}+z^{-2i})\sum_{j=i}^\infty (\mp 1)^{j-i} \binom{2j}{j-i} V^{\pm}_{j,m}(q) \\
		& \qquad = \sum_{i=-\infty}^\infty z^{2i} \sum_{j=|i|}^\infty (\mp 1)^{j-|i|} \binom{2j}{j-|i|} V^{\pm}_{j,m}(q).
	\end{align*}
	Thus we deduce that
	\begin{align*}
		& \sum_{i=-\infty}^\infty z^{2i} \sum_{j=|i|}^\infty (\mp 1)^{j-|i|} \binom{2j}{j-|i|} V^{\pm}_{j,m}(q)\\
		& \qquad =(\pm q;q)^2_m \sum_{i=-\infty}^\infty z^{2i} \sum_{j=0}^\infty 
		{m-1+j \brack j} {m-1+|i|+j \brack |i|+j} q^{|i|+2j}.
	\end{align*}
	By the equating the coefficient of $z^{2k}$ ($k\geq 0$) in this identity, we obtain
	\begin{align*}
		& \sum_{j=k}^\infty (\mp 1)^{j-k} \binom{2j}{j-k} V^{\pm}_{j,m}(q)\\
		&\qquad =(\pm q;q)^2_m \sum_{j=0}^\infty 
		{m-1+j \brack j} {m-1+k+j \brack k+j} q^{k+2j}.
	\end{align*}    
	This concludes the proof of the first identity.

    The second identity follows in a similar manner.
	So considering \eqref{Cauchy}, we can write:
	\begin{align*}
	& \sum_{i=0}^m W^{\pm}_{i,m}(q)\, (z\mp z^{-1})^{2i} \\ 
	&\qquad   = \prod_{i=1}^m \left(1-\frac{q^{2i-1}}{(1\mp q^{2i-1})^2} (z\mp z^{-1})^2 \right)^{-1}\\
	&\qquad   = (\pm q;q^2)^2_m \prod_{i=1}^m \frac{1}{(1-q^{2i-1}z^2)(1-q^{2i-1}z^{-2})} \\
	&\qquad   = \frac{(\pm q;q^2)^2_m}{(qz^2;q^2)_m\,(qz^{-2};q^2)_m} \\
	%&\qquad   = \frac{(qz^{-2};q)_m}{(\pm q;q)^2_m} \sum_{i=0}^m (-1)^i q^{\binom{i+1}{2}} {m \brack i} z^{2i} \\
	%\intertext{\hfill by \eqref{Rbt} with $t$ replaced by $-qz^2$}
	&\qquad   = (\pm q;q^2)^2_m \left(\sum_{i=0}^m {m-1+i \brack i}_{q^2} q^i\,z^{2i} \right)
	\left(\sum_{i=0}^m {m-1+i \brack i}_{q^2} q^i\,z^{-2i} \right)\\
	\intertext{\hfill by \eqref{Cauchy} with $q$ replaced by $q^2$ and $t$ replaced by $-qz^{\pm2}$}
%	&\qquad   = \frac{1}{z^{2m}\,(\pm q;q^2)^2_m} \left(\sum_{i=0}^m (-1)^i\, q^{i^2} {m \brack i}_{q^2} z^{2i} \right)
%	\left(\sum_{i=0}^m (-1)^{m-i}\, q^{(m-i)^2} {m \brack m-i}_{q^2} z^{2i} \right)\\
%	&\qquad   = \frac{1}{z^{2m}\,(\pm q;q^2)^2_m} \sum_{i=0}^{2m} \sum_{j=0}^i (-1)^{m-i}\, q^{j^2+(m-i+j)^2} {m\brack j}_{q^2} {m\brack m-i+j}_{q^2} z^{2i}\\
%	&\qquad   = \frac{1}{z^{2m}\,(\pm q;q^2)^2_m} \sum_{i=0}^{2m} \sum_{j=0}^i (-1)^{m-i}\, q^{(i-j)^2+(m-j)^2} {m\brack i-j}_{q^2} {m\brack m-j}_{q^2} z^{2i}.\\
		&\qquad = (\pm q;q^2)^2_m \sum_{i=-\infty}^\infty z^{2i} \sum_{j=0}^\infty 
	{m-1+j \brack j}_{q^2} {m-1+|i|+j \brack |i|+j}_{q^2} q^{|i|+2j}.
	\end{align*}
	On the other hand, we have:
	\begin{align*}
& \sum_{i=0}^\infty W^{\pm}_{i,m}(q)\, (z\mp z^{-1})^{2i}  = \sum_{i=-\infty}^\infty z^{2i} \sum_{j=|i|}^\infty (\mp 1)^{j-|i|} \binom{2j}{j-|i|} W^{\pm}_{j,m}(q).
\end{align*}
	Thus we deduce that
	\begin{align*}
	& \sum_{i=-\infty}^\infty z^{2i} \sum_{j=|i|}^\infty (\mp 1)^{j-|i|} \binom{2j}{j-|i|} W^{\pm}_{j,m}(q)\\
	& \qquad = (\pm q;q^2)^2_m \sum_{i=-\infty}^\infty z^{2i} \sum_{j=0}^\infty 
	{m-1+j \brack j}_{q^2} {m-1+|i|+j \brack |i|+j}_{q^2} q^{|i|+2j}.
	\end{align*}
	By equating the coefficient of $z^{2k}$ ($k\geq 0$) in this identity, we obtain
	\begin{align*}
	& \sum_{j=k}^\infty (\mp 1)^{j-k} \binom{2j}{j-k} W^{\pm}_{j,m}(q)\\
	& \qquad = (\pm q;q^2)^2_m \sum_{j=0}^\infty 
	{m-1+j \brack j}_{q^2} {m-1+k+j \brack k+j}_{q^2} q^{k+2j}.
	\end{align*}
	This completes the proof of the second identity.
}

\section{Proof of Theorem \ref{TT2}}

Let sequences $(V_j)_{j\geq 0}$ and $(R_k)_{k\geq 0}$ be related by
\begin{equation}\label{eq:relation}
	R_k=\sum_{j=k}^{\infty}(-1)^{\,j-k}\binom{2j}{\,j-k\,}\,V_j,\qquad k\geq 0,
\end{equation}
(with the sums convergent, e.g. finite-support $V_j$). Define ordinary generating functions
\[
V(x)=\sum_{j=0}^\infty V_j\, x^j,\qquad R(y)=\sum_{k=0}^\infty R_k\, y^k.
\]
We can write:
\begin{align*}
	R(y) &= \sum_{k=0}^\infty R_k\,y^k \\
	&= \sum_{k=0}^\infty \sum_{j=k}^\infty (-1)^{j-k}\,\binom{2j}{j-k}V_j\,y^k \\
	&= \sum_{k=0}^\infty \sum_{j=0}^\infty (-1)^{j}\,\binom{2j+2k}{j}V_{j+k}\,y^k \\
	&= \sum_{j=0}^\infty V_j \sum_{k=0}^j (-1)^{j-k}\,\binom{2j}{j-k}\,y^k \\
	&= \sum_{j=0}^\infty V_j\,y^j \sum_{r=0}^{2j} \binom{2j}{r}\left(-\frac{1}{y}\right)^{r} \\
	&= \sum_{j=0}^\infty V_j\,\frac{(y-1)^{2j}}{y^j} \\
	&= V\big(\phi(y)\big),
\end{align*}
where
$$\phi(y)=\frac{(y-1)^2}{y}.$$
%Our task is to invert the map \(V\mapsto R\) given by composition with \(\phi\).

We recover the coefficients \(V_j\) from \(R\) by Cauchy's coefficient formula \cite[Theorem 23]{Kaplan2002}
applied to the composite \(V(x)=R(\Phi(x))\), where \(\Phi\) denotes a local
(analytical) inverse of \(\phi\) near the point \(y=1\). Equivalently, 
without explicitly solving for \(\Phi\), 
we use the change of variables \(x=\phi(y)\)
in the Cauchy integral for the coefficients of \(V\):
\[
\begin{aligned}
	V_j &= [x^j]V(x)
	= \frac{1}{2\pi i}\oint_{|x|=\epsilon} \frac{V(x)}{x^{j+1}}\,dx
	= \frac{1}{2\pi i}\oint_{y=1} \frac{R(y)}{\phi(y)^{\,j+1}}\,\phi'(y)\,dy.
\end{aligned}
\]
Substitute \[R(y)=\sum_{k=0}^\infty R_k\, y^k,\] interchange sum and integral (finite
support justifies this), and extract the contribution of \(R_k\):
$$
V_j = \sum_{k=0}^\infty B_{k,j}\,R_k ,
$$
where the inverse kernel is given by the contour integral
\[
B_{k,j} = \frac{1}{2\pi i}\oint_{y=1} \frac{y^k}{\phi(y)^{\,j+1}}\,\phi'(y)\,dy.
\]
Since \(\phi(y)=(y-1)^2/y\) we get
\[
\phi'(y)=\frac{(y-1)(y+1)}{y^2}.
\]
Thus the integrand simplifies as
\[
\frac{y^k}{\phi(y)^{\,j+1}}\,\phi'(y)
= y^{k+j-1}\,(y+1)\,(y-1)^{-2j-1}.
\]
Change variables to \(z=y-1\). Then \(y=1+z\) and the integration contour is
a small circle around \(z=0\). The integrand becomes
\[
(1+z)^{\,j+k-1}\,(2+z)\,z^{-2j-1}.
\]
Hence by the residue theorem the integral picks out the coefficient of
\(z^{2j}\) in the (polynomial) expansion of \((2+z)(1+z)^{\,j+k-1}\). That is,
\[
 \; B_{k,j} = [z^{2j}]\,(2+z)\,(1+z)^{\,j+k-1}. \;
\]
Expanding the binomial series gives for \(j\geq 1\):
\[
\begin{aligned}
	B_{k,j}
	&= 2\binom{j+k-1}{2j} + \binom{j+k-1}{2j-1}\\
	&= \frac{2k}{k+j}\binom{k+j}{2j},
\end{aligned}
\]
while the exceptional row \(j=0\) is obtained by direct inspection:
\(B_{0,0}=1\) and \(B_{k,0}=2\) for \(k\geq 1\) (indeed
\([z^0](2+z)(1+z)^{k-1}=2\) for \(k\geq 1\)).
By \eqref{eq:relation}, with $R_k$ replaced by $(-1)^k\,R'_k$ and $V_j$ replaced by $(-1)^j\,V'_j$, we obtain
$$R'_k=\sum_{j=k}^{\infty} \binom{2j}{\,j-k\,}\,V'_j,\qquad k\geq 0,$$
It is clear that
$$V'_j = \sum_{k=j}^\infty (-1)^{k-j}\,B_{k,j}\,R'_k,\qquad j\geq 0.$$
These yield the explicit formulas for all \(j,k\) with \(k\geq j\) stated in
the theorem.

\section{Proof of Theorem \ref{T4}}

The limiting case $m\to\infty$ of Theorem \ref{T1} reads as follows: 
	\begin{align*}
& \sum_{j=k}^\infty (\mp 1)^{j-k} \binom{2j}{j-k} V^{\pm}_{j}(q)
=(\pm q;q)^2_\infty \sum_{j=0}^\infty 
\frac{q^{2j+k}}{(q;q)_j\,(q;q)_{j+k}}
\end{align*} 
and
\begin{align*}
& \sum_{j=k}^\infty (\mp 1)^{j-k} \binom{2j}{j-k} W^{\pm}_{j}(q)
=(\pm q;q^2)^2_\infty \sum_{j=0}^\infty 
\frac{q^{2j+k}}{(q^2;q^2)_j\,(q^2;q^2)_{j+k}}.
\end{align*} 

The first identity of Theorem \ref{T4} follows easily if we consider the following result for wich we present a detailed proof, with all algebraic reindexings made explicit.

\medskip
\begin{lemma}\label{Lemma1}
	Let $k$ be non-negative integer. For $|q|<1$,
%	\begin{align*}
%	& (q;q)^2_\infty \sum_{j=0}^\infty 
%	\frac{q^{2j+k}}{(q;q)_j\,(q;q)_{j+k}} = -1+\frac{(-1)^k\,(1+q^k)}{q^{k(k+1)/2}} \sum_{j=k}^\infty (-1)^j\,q^{j(j+1)/2}.
%	\end{align*} 
	\begin{align*}
	& (q;q)^2_\infty \sum_{j=0}^\infty 
	\frac{q^{2j+k}}{(q;q)_j\,(q;q)_{j+k}} = -1+(1+q^k) \sum_{j=k}^\infty (-1)^{j-k}\,q^{j(j+1)/2-k(k+1)/2}.
	\end{align*} 
\end{lemma}
\medskip

\begin{proof} We shall make use of Euler’s expansion of the infinite $q$-Pochhammer product \cite[Theorem~27]{Johnson}:
	for $|q|<1$,
	\begin{equation} \label{Euler1}
	(z;q)_\infty
	= \sum_{j=0}^\infty 
	\frac{(-1)^j\, q^{j(j-1)/2}}{(q;q)_j}\, z^{j}.
	\end{equation}
	We also recall the following theorem of Euler \cite[Theorem~25]{Johnson}: if $|q|<1$ and $z<1$, then
	\begin{equation} \label{Euler2}
	\sum_{j=0}^\infty \frac{z^{j}}{(q;q)_j}
	= \frac{1}{(z;q)_\infty}.
	\end{equation}
	On the other hand, for integers $j,k\geq 0$ we note that
	\[
	(q;q)_{j+k}
	= \frac{(q;q)_\infty}{(q^{\,j+k+1};q)_\infty},
	\qquad\text{and hence}\qquad
	\frac{1}{(q;q)_{j+k}}
	= \frac{(q^{\,j+k+1};q)_\infty}{(q;q)_\infty}.
	\]
	So the left-hand side of our identity can be written as:
	\begin{align*}
		LHS
		&= (q;q)_\infty^2 \sum_{j=0}^\infty \frac{q^{2j+k}}{(q;q)_j (q;q)_{j+k}} \\
		&= (q;q)_\infty \sum_{j=0}^\infty \frac{q^{2j+k}}{(q;q)_j} (q^{\,j+k+1};q)_\infty\\
		&= (q;q)_\infty \sum_{j=0}^\infty \frac{q^{2j+k}}{(q;q)_j}
		\sum_{n=0}^\infty \frac{(-1)^n\, q^{n(n-1)/2 + n(j+k+1)}}{(q;q)_n} \\
		\intertext{\hfill by \eqref{Euler1} with $z$ replaced by $q^{j+k+1}$}
		&= (q;q)_\infty \sum_{n=0}^\infty \frac{(-1)^n\, q^{n(n-1)/2 + n(k+1)}}{(q;q)_n}
		\sum_{j=0}^\infty \frac{q^{j(2+n)+k}}{(q;q)_j}\\
		&= (q;q)_\infty \sum_{n=0}^\infty \frac{(-1)^n\, q^{n(n-1)/2 + n(k+1)}\, q^k}{(q;q)_n\, (q^{n+2};q)_\infty}\\
			\intertext{\hfill by \eqref{Euler2} with $z$ replaced by $q^{2+n}$}
		&= q^k \sum_{n=0}^\infty (-1)^n\, q^{n(n-1)/2 + n(k+1)}\, \frac{(q;q)_{n+1}}{(q;q)_n} \\
		&= q^k \sum_{n=0}^\infty (-1)^n\, q^{n(n+2k+1)/2}\, (1-q^{n+1}) \\
		&= q^k \sum_{n=0}^\infty (-1)^n\, q^{n(n+2k+1)/2}
		- q^k \sum_{n=0}^\infty (-1)^n\, q^{n(n+2k+1)/2 + n + 1}\\
		&= q^k \sum_{n=0}^\infty (-1)^n\, q^{n(n+2k+1)/2}
		- q^k \sum_{n=1}^\infty (-1)^{n-1} q^{n(n+2k+1)/2-k}\\
		&= q^k + q^k\sum_{n=1}^\infty (-1)^n\, q^{n(n+2k+1)/2}
- \sum_{n=1}^\infty (-1)^{n-1} q^{n(n+2k+1)/2}\\
		&= q^k + (1+q^k)\sum_{n=1}^\infty (-1)^n\, q^{n(n+2k+1)/2}\\
		&= q^k + \frac{1+q^k}{q^{k(k+1)/2}}\sum_{n=1}^\infty (-1)^n\, q^{(n+k)(n+k+1)/2}\\
		&= q^k + \frac{(-1)^k\,(1+q^k)}{q^{k(k+1)/2}}\sum_{n=k+1}^\infty (-1)^n\, q^{n(n+1)/2}\\
		&= q^k + \frac{(-1)^k\,(1+q^k)}{q^{k(k+1)/2}}\left(-(-1)^k\,q^{k(k+1)/2} + \sum_{n=k}^\infty  (-1)^n\, q^{n(n+1)/2} \right)\\
		&= -1 + \frac{(-1)^k\,(1+q^k)}{q^{k(k+1)/2}} \sum_{n=k}^\infty  (-1)^n\, q^{n(n+1)/2}.
	\end{align*}
This completes the proof of the identity.
\end{proof}

The second identity of Theorem \ref{T4} follows easily if we consider the following result for wich we present a detailed proof, with all algebraic reindexings made explicit.

\medskip
\begin{lemma}\label{Lemma2}
	Let $k$ be non-negative integer. For $|q|<1$,
\begin{align*}
& (q^2;q^2)^2_\infty \sum_{j=0}^\infty 
\frac{q^{2j+k}}{(q^2;q^2)_j\,(q^2;q^2)_{j+k}} = \sum_{j=k}^\infty (-1)^{j-k}\,q^{j(j+1)-k^2}.
\end{align*} 
\end{lemma}
\medskip

\begin{proof} The proof of this lemma is quite similar to the proof of Lemma \ref{Lemma1}. So we can write:
\begin{align*}
& (q^2;q^2)^2_\infty \sum_{j=0}^\infty 
\frac{q^{2j+k}}{(q^2;q^2)_j\,(q^2;q^2)_{j+k}} \\
& \qquad = (q^2;q^2)_\infty \sum_{j=0}^\infty 
\frac{q^{2j+k}}{(q^2;q^2)_j} \,(q^{2(j+k)+2};q^2)_\infty \\
& \qquad = (q^2;q^2)_\infty \sum_{j=0}^\infty \frac{q^{2j+k}}{(q^2;q^2)_j}
\sum_{n=0}^\infty \frac{(-1)^n\, q^{n(n-1)}\, q^{2n(j+k+1)}}{(q^2;q^2)_n} \\
\intertext{\hfill by \eqref{Euler1} with $z$ replaced by $q^{2(j+k)+2}$}
& \qquad = (q^2;q^2)_\infty \sum_{n=0}^\infty \frac{(-1)^n\, q^{n(n-1)+2n(k+1)}}{(q^2;q^2)_n}
\sum_{j=0}^\infty \frac{q^{j(2+2n)+k}}{(q^2;q^2)_j} \\
& \qquad = (q^2;q^2)_\infty \sum_{n=0}^\infty \frac{(-1)^n\, q^{n(n-1)+2n(k+1)}\,q^k}{(q^2;q^2)_n\,(q^{2n+2};q^2)_{\infty}}\\
\intertext{\hfill by \eqref{Euler2} with $z$ replaced by $q^{2n+2}$}
& \qquad = q^k \sum_{n=0}^\infty (-1)^n\, q^{n(n-1)+2n(k+1)} \\
& \qquad = \sum_{n=k}^\infty (-1)^{n-k}\, q^{n(n-1)-k^2}.
\end{align*}	
This completes the proof of the identity.	
\end{proof}

\section{Concluding remarks}

In this work, we investigated the MacMahon-type $q$-series $V_k^{\pm}(q)$ and $W_k^{\pm}(q)$, together with their truncated forms $V_{k,m}^{\pm}(q)$ and $W_{k,m}^{\pm}(q)$. Through a combination of generating function manipulations, Gaussian polynomial identities, and basic hypergeometric series representations, we established new structural relations satisfied by these series. In particular, Theorems \ref{T1} and \ref{TT2} express linear combinations of $V_{k,m}^{\pm}(q)$ and $W_{k,m}^{\pm}(q)$ explicitly in terms of $q$-binomial coefficients, thus revealing refined analytic behavior of these functions.

Passing to the limit $m \to \infty$, Theorem \ref{T4} yields new product-false theta decompositions for $V_k^{\pm}(q)$ and $W_k^{\pm}(q)$, connecting their coefficients to well-studied partition families, including overpartition pairs and bipartitions with odd parts distinct. These identities also imply new positivity properties for weighted convolution sums involving $v_k^{\pm}(n)$ and $w_k^{\pm}(n)$.

Finally, Theorem \ref{TT4} provides explicit expansions of 
$$V_j^{\pm}(q)=(\pm q;q)^2_\infty \sum_{k=j}^\infty \sum_{i=0}^\infty \frac{(\pm 1)^{k-j}\,B_{k,j}\,q^{2i+k}}{(q;q)_i\,(q;q)_{i+k}}$$ 
and 
$$W_j^{\pm}(q)=(\pm q;q^2)^2_\infty \sum_{k=j}^\infty \sum_{i=0}^\infty \frac{(\pm 1)^{k-j}\,B_{k,j}\,q^{2i+k}}{(q^2;q^2)_i\,(q^2;q^2)_{i+k}}$$ 
in terms of alternating sums of triangular-number exponents, leading to new representations for squares of Ramanujan theta functions. In particular, the case $j=0$ yields weighted enumerations for ordered representations of integers as sums of two squares or two triangular numbers, revealing unexpected connections between these functions and classical representation theory. Similar results can be obtained for the sum–of–divisors functions by considering the case 
$j=1$. For instance, the generating function for the sum of the divisors of $n$ admits a new representation involving the squares of natural numbers:
$$\sum_{k=1}^\infty \frac{k\,q^k}{1-q^k}=(q;q)^2_\infty \sum_{k=1}^\infty \sum_{i=0}^\infty \frac{k^2\,q^{2i+k}}{(q;q)_i\,(q;q)_{i+k}}.$$ 

The results presented here suggest several promising avenues for further study. These include the exploration of higher-rank analogues of our constructions, modular and quasimodular properties of the resulting series, and arithmetic questions such as congruences and parity distributions of the associated partition functions. We hope that the structural discoveries developed in this paper will inspire continued investigation into MacMahon-type $q$-series and their deep connections with partitions, hypergeometric functions, and the theory of theta functions.

\section*{Declaration of competing interest}

The author declares that he has no known competing financial interests or personal relationships that could appear to influence the work reported in this paper.

\section*{Data availability}

No data was used for the research described in the article.



%%       \citet{<label>} ==> Jones et al. [21]
%%       \citep{<label>} ==> [21]
%%

  \bibliographystyle{elsarticle-num-names} 
  \bibliography{bibliography}


%\begin{thebibliography}{00}



%\end{thebibliography}
\end{document}